\documentclass[11pt]{article}
\usepackage[cp866]{inputenc} 
\usepackage[russian]{babel}

\begin{document} 

\centerline{\Large\bf Quantum Properties of } 
\centerline{\Large\bf Mathematical Physics Equations.}
\centerline{\Large\bf Generation of Quantum Structures }
\centerline {\bf L.I. Petrova}
\centerline{\it Lomonosov Moscow State University, Faculty of }
\centerline{\it Computational Mathematics and Cybernetics, e-mail: ptr@cs.msu.ru}

\renewcommand{\abstractname}{Abstract} 
\begin{abstract} 

It is shown with the help of skew-symmetric forms that the mathematical physics equations, on which no additional conditions are imposed, have quantum properties. And this is due to the integrability properties of differential equations, which depends on the consistency of derivatives with respect to different variables and the consistency of equations, if the mathematical physics equations are a system of equations.
It was found that such equations on the original tangent space turn out to be non-integrable. Their derivatives do not form a differential. The integrability of such equations is realized only on the structures of a cotangent integrable manifold. This happens using a degenerate, non-differential-preserving transformation that has quantum properties. When implementing degenerate transformations, mini structures (quanta) arise, from which integrable structures are formed.
Such properties of integrability of mathematical physics equations and features of degenerate transformations reveal the quantum properties of mathematical physics equations and their ability to generate quantum structures.

\end{abstract} 


\bigskip

{\Large\bf Introduction} 

The theory of skew-symmetric differential forms made it possible to reveal the unique capabilities of differential mathematical physics equations, on which no additional conditions are imposed. Such equations have quantum properties and can generate quantum structures.

\bigskip 

Quantum properties of the mathematical physics equations are related to the integrability of differential equations. This was attained by studying the integrability of the mathematical physics equations, which, as shown, depends on the consistency of the differential equations derivatives with respect to different variables and on the consistency of the equations if the mathematical physics equations are a system of equations. (Derivatives are consistent if they form a differential.) 

Usually, when solving mathematical physics equations, issues of consistency of differential equations derivatives and consistency of equations forming a system of equations are not taken into account. But it turns out that this is what integrability of differential equations and their quantum properties depend on.

The study of the consistency of derivatives and the study of the consistency of equations forming a system of equations showed that the integrability of equations has specific features. 

On tangent space, the mathematical physics equations turn out to be non-integrable. Derivatives of equations do not form a differential. In this case, solutions to equations in coordinate space are not functions.
The mathematical physics equations can become integrable only if additional conditions are being implemented. 

 And here there is a very specific detail. 
  
The integrability of such equations is realized not on the tangent space, since it is non-integrable, but on the cotangent, integrable manifold. In this case, an integrable structure of the cotangent manifold forms and a transition occurs from the tangent space to the structures of the cotangent manifold. On the emerging integrable structure, solutions to the mathematical physics equations become discrete functions. 

It turns out that the mathematical physics equations have solutions that are defined on different spatial objects: on the original coordinate space and on the cotangent manifold structures. 

How is the connection between such dual solutions is realized?

It turns out that there are mathematical operators that describe the relationship between solutions on different spatial objects. 

These are degenerate, non-differential-preserving transformations of skew-symmetric forms: transformations in which there is a transition from a differential that is not equal to zero to the differential that equals zero.

The degenerate transformation has a unique feature that is not found in other mathematical formalisms. When implementing a degenerate transformation, an integrable structure arises and a transition occurs from a non-integrable space to the structures of an integrable manifold. This makes it possible to describe discrete transitions and the emergence of various structures.

Degenerate transformations occur discretely when implementing any degrees of freedom. This corresponds to the vanishing of such functional expressions as Jacobians, determinants, Poisson brackets, residues, etc. Such conditions describe the realized integrable structures. 

The degenerate transformation brings the mathematical physics equations to an integrable form. 

An example of a degenerate transformation is the Legendre transformation. As is known, with the Legendre transformation, Hamiltonian structures emerge and a transition occurs from a Lagrangian, non-integrable, manifold to a Hamiltonian, integrable, manifold. This brings Euler's equation to an integrable form. The Legendre transformation was also used to reduce differential equations to integrable form in the works of Hilbert and Courant and to obtain the integrable Hamilton-Jacobi equation. 

But at the same time, the nonintegrability of the mathematical physics equations was usually not studied. 

Namely, the solutions that are obtained when differential equations are not integrable are the generator of quantum structures that arise when the integrability of differential equations is realized. And this is described by a degenerate transformation. 

As will be shown in this work, degenerate transformations reveal the processes of emergence of various structures and invariant objects. 

\bigskip

Quantum properties of differential mathematical physics equations, which are due to the integrability of differential equations and are associated with the properties of degenerate transformations, were shown using the example of a first-order partial differential equation (Chapter 1) and using the example of mathematical physics equations that consist of several equations, that is are a system of equations (Chapter 2). Some principles for the implementation of quantum processes that reveal the mechanism of origin of existing objects and their quantum structure are described in Chapter 3.

\section {Dependence of the integrability of differential equations on the consistency of derivatives. Discrete properties of differential mathematical physics equations} 

This paper examines differential mathematical physics equations that describe physical processes and on which no additional conditions are imposed. 

The integrability of such differential equations, as the study shows, depends on the consistency of the derivatives of the differential equations. (Derivatives are consistent if they form a differential.) Moreover, if the mathematical physics equations are a system of equations, then the integrability of such a system depends on the consistency of the equations forming this system. 

\bigskip

Here it should be emphasized that the study of integrability of differential equations is based on the properties of skew-symmetric differential forms  [1-3]. It is obvious that the characteristic properties of integrable equations are the presence of differentials and the ability to reduce the equation to an identical relation in differentials, which will allow one to explicitly integrate the differential equation. Skew-symmetric differential forms have these properties. Unlike other existing mathematical formalisms, they can take the form of differentials and differential expressions. Appendix 1 shows some properties of skew-symmetric forms that will be used in this study. 

\bigskip

This chapter will show the dependence of the integrability of differential equations on the consistency between the derivatives of differential equations, and the related features of solutions of differential equations.

\subsection*{Study of consistency of differential equations derivatives.}

The consistency between the derivatives of the sought functions in differential equations can be seen using the example of a first-order partial differential equation:
$$ F(x^i,\,u,\,p_i)=0,\quad p_i\,=\,\partial u/\partial x^i \eqno(1)$$ 

The derivatives of a differential equation are consistent if they form a differential  $ du\,$. That is, if the relation 
$ du\, = \,p_i\,dx^i$ is satisfied (here we mean summation over repeating indices). 

If the differential formed by the derivatives $p_i\,$ is written in the form of a skew-symmetric form 
$\theta\,=\,p_i\,dx^i$, then the condition for consistency of derivatives can be written in the form
$$ du\,=\,\theta\eqno(2)$$. 

However, in the general case it does not follow from the differential equation that the derivatives 
$p_i\,=\,\partial u/\partial x^i $, satisfying equation (1), form a differential. That is, it does not follow from equation (1). that the skew-symmetric form  $\theta\,=\,p_i\,dx^i$ must be a differential and satisfy relation (2).

A skew-symmetric form can be a differential if it is a closed form, that is, the condition for the form to be closed is satisfied: the differential of the form is equal to zero and the differential of the corresponding dual form is equal to zero if the form is a closed inexact form. (Appendix 1, [1]). 

Form $\theta\,=\,p_i\,dx^i$ can be a differential only if the conditions of closedness of the form 
$\theta\,=\,p_idx^i$ and the corresponding dual form are satisfied (in this case, the functional $F$ plays the role of the dual form for $\theta $)): that is, the vanishing of the differentials of the skew-symmetric form and the corresponding dual form:
$$\cases {dF(x^i,\,u,\,p_i)\,=\,0\cr
d(p_i\,dx^i)\,=\,0\cr}\eqno(3)$$

If we expand the differentials, we get a system of homogeneous equations for $ dx ^ i $ ahd $ dp_i $  (that is, in cotangent space):
$$\cases {\displaystyle \left ({{\partial F}\over {\partial x^i}}\,+\,
{{\partial F}\over {\partial u}}\,p_i\right )\,dx^i\,+\,
{{\partial F}\over {\partial p_i}}\,dp_i \,=\,0\cr
dp_i\,dx^i\,-\,dx^i\,dp_i\,=\,0\cr} \eqno(4)$$

The condition for the solvability of this system, that is {\bf vanishing of the determinant composed of the coefficients of $dx^i$, $dp_i$}, has the following form:
$$
{{dx^i}\over {\partial F/\partial p_i}}\,=\,{{-dp_i}\over 
{\partial F/\partial x^i+p_i\partial F/\partial u}} \eqno (5)
$$ 

This condition defines an integrable structure: a pseudostructure on which relations (3) are satisfied. Such an integrable structure is defined on a cotangent manifold. (An example of such an integrable structure is the characteristics of differential equations.) 

The fulfillment of relation (3) indicates that on the integrable structure the form $\theta \,=\,p_i\,dx^i$, formed by the derivatives of the differential equation, turns out to be a closed form, i.e. {\bf becomes a differential}. The derivatives of the differential equation turn out to be consistent. So, relation (2) is satisfied.

It turns out that the derivatives of the differential equation (1) turn out to be consistent only on the structures of the cotangent manifold. 

And on the tangent space of the differential equation (1) the derivatives turn out to be inconsistent. This indicates that the tangent space turns out to be non-integrable.

\subsection*{Dependence of the integrability of differential equations on the consistency of derivatives.} 

It was shown that on integrable structures of a cotangent manifold, the derivatives of a differential equation turn out to be consistent, so, they form a differential. 

In this case, relation (2) becomes an identical relation. Skew-symmetric form $\theta\,=\,p_i\,dx^i$, standing on the right side of relation (2), is a differential, just like the left side of the relation. Since relation (2) consists only of a differential, it can be integrated and obtain a solution  $u$ to the differential equation (1),  which is a function, since the derivatives of this solution form a differential. 

Such solutions are discrete functions because they are defined on integrable structures. 

It turns out that equation (1), if the derivatives are consistent, becomes integrable: its solutions become functions. {\bf And this is realized on integrable structures of a cotangent manifold}. It should be emphasized that the consistency of the derivatives of the differential equation is associated with the fulfillment of relation (3), which is a condition for the closedness of the form $\theta \,=\,p_i\,dx^i$ and occurs when an additional condition is realized: when the determinant of the system (4) vanishes, which determines the integrable structure (5) of the cotangent manifold. 

{\bf And on the original tangent space, differential equation (1) turns out to be non-integrable.} 

In this case, relation (3), which is the condition for the closedness of the form $\theta \,=\,p_i\,dx^i$, is not true. 

This means that the form $\theta \,=\,p_i\,dx^i$, formed by the derivatives of a differential equation, is not a closed form, so, it is not a differential. The derivatives of equation (1) turn out to be inconsistent. 

In this case, relation (2) turns out to be non-identical, since on the left there is a differential, and on the right, there is a skew-symmetric form, which is not a differential. It is impossible to directly integrate such a relation  and obtain a solution, which is a function. That is, equation (1) on the tangent space turns out to be non-integrable. The solution to equation (1) on the original coordinate space is not a function: the derivatives of such a solution do not form a differential. But it must be emphasized that such a solution has a physical meaning, since the equation under study was assumed to describe real physical processes. (Such solutions can be obtained by numerical methods.) 

(Here it should be noted that the non-integrability of differential equations on tangent space is due to the fact that they include terms associated with non-potential influences to which the described system is subjected. This also explains the fact that the solutions of non-integrable equations are not functions.) 

\bigskip

It turns out that differential equation (1) has solutions of different types: solutions on coordinate space, which are not functions, and solutions on integrable structures of the cotangent manifold, which are discrete functions. 

And this, as shown, is conected to the integrability of the differential equation. 

It can be noted that such properties are possessed by differential equations that describe physical processes and on which integrability conditions are not imposed. Such conditions, as will be shown, are realized in the process of solving equations. 

\subsection*{Discrete properties of differential mathematical physics equations.} 

Here you should pay attention to the peculiarity of double solutions. They are implemented on different spatial objects. And this, as shown, is related to the integrability of the differential equation. 

The question arises. How does the transition from a solution of the first type to a solution of the second type occur? That is, how does the transition from tangent space to integrable structures of a cotangent manifold occur? 

It was shown that the type of solutions depends on whether the derivatives of the differential equation form a differential or not. So, whether relation (3) is satisfied or not. 

As has been shown, solutions in coordinate space are not functions, because the derivatives of the differential equation (1) are inconsistent. This is due to the fact that relation (3) does not hold: differential $d(p_i\,dx^i)\,$ is not equal to zero. 

And solutions on integrable structures turn out to be functions (discrete), because relation (3) is satisfied: differential $d(p_i\,dx^i)\,$ is equal to zero. 

It turns out that the transition from a solution in coordinate space to a solution in integrable structures occurs {\bf as a transition from a differential not equal to zero to a differential equal to zero}. 

That is, it happens {\bf degenerate transformation: a transformation that does not preserve differential}. 

A degenerate transformation has a unique property. Under a degenerate transformation, an integrable structure arises and a transition occurs from a non-integrable space to the structures of an integrable manifold.

{\footnotesize[An example of degenerate transformations, as noted, is the Legendre transformation. As is known, the Legendre transformation leads to the emergence of Hamiltonian structures and the transition from a non-integrable Lagrangian manifold to an integrable Hamiltonian manifold [4]. 
As could be seen when studying equation (1), the degenerate transformation leads the differential equation to an integrable form. This property of the degenerate Legendre transform, as already noted, was used in the work of Courant and Hilbert [5] when calculating the characteristics of a differential equation on which the integrability of the equation is realized. The integrability of the Hamilton-Jacobi equation is based on the Legendre transform.]}
 
Degenerate transformations occur when any degrees of freedom are realized, that is, when additional conditions are realized, which are integrability conditions. Such conditions determine the emerging integrable structure.

For equation (1), as was shown, such an additional condition is {\bf the vanishing of the determinant of system (4)}  (it determines the integrable structure (5) of the cotangent manifold).

In the general case, these conditions are associated with the vanishing of such functional expressions as determinants, Jacobians, Poisson brackets, and residues. 

It turns out that degenerate transformations, in which the mathematical physics equations become integrable, can only occur discretely when implementing any degrees of freedom. This occurs during the solution process and is accompanied by the emergence of integrable structures with solutions that are discrete functions. 

Integrable structures, in which the solutions to equations are discrete functions, have a physical meaning. These are physical structures that describe invariant objects. 

Such properties of degenerate transformations reveal discrete properties of differential equations. 

In the next chapter, when studying the mathematical physics equations, which consist of several equations, that is, they are a system of equations, the process of implementing degenerate transformations will be described in detail, and the unique role of degenerate transformations for the mathematical physics equations, revealing their quantum properties, will be shown.

\section {The dependence of the integrability of the mathematical physics equations, which are a system of equations, on the consistency of the equations included in the system. Discrete properties of mathematical physics equations} 

This chapter examines the differential mathematical physics equations, which are a system of equations. It is shown that the integrability of such a system depends on the consistency of the equations forming this system.

Here it should be emphasized that equations are being studied that are not subject to additional conditions (such conditions, as shown, are spontaneously realized in the process of solving mathematical physics equations). Equations on which no additional conditions are imposed are presented, for example, in [6-9]. (In physics, equations with integrability conditions are practically used. A broad overview of such works is given in the book [10].)

This study was carried out using the example of equations of mechanics and physics of continuous media, such as thermodynamic, gas-dynamic, cosmological, electromagnetic, etc. Such mathematical physics equations consist of equations of the laws of conservation of energy, momentum, angular momentum and mass.

\subsection*{Study of the consistency of conservation laws equations.} 

Let us consider the consistency of the equations of the laws of conservation of energy and momentum, which are included in the mathematical physics equations. 

It should be noted here that the present study has its own peculiarities. The equations of conservation laws are transformed into equations for functionals that characterize the state of the described medium. Functionals such as wave function, entropy, action functional, etc. (which are also field theory functionals) are such state functionals [11]. 

In a coordinate system associated with a manifold formed by particle trajectories (an example of such a coordinate system is the Lagrangian coordinate system), the energy and momentum equations can be written as:
$$ 
{{\partial \psi }\over {\partial \xi ^1}}\,=\,A_1 \eqno(6)
$$ 
$$
{{\partial \psi}\over {\partial \xi^{\nu }}}\,=\,A_{\nu },\quad \nu \,=\,2,\,...\eqno(7)
$$ 
Here  $\psi$ is the state functional,  $\xi^1$ and $\xi ^{\nu }$ are coordinates along the trajectory and in the direction normal to the trajectory, $A_1$ is a quantity depending on the characteristics of the material environment being described and on external energy influences on the environment, and  $A_{\nu }$ are quantities depending on the characteristics of the material environment and external force influences. 

{\footnotesize [So the energy equation, expressed through the action functional $S$, has the same form: 
${{\partial \,S }/{\partial \xi ^1}}\,=\,L $, where  $\psi \,=\,S$, $A_1\,=\,L$ is the Lagrange function [4]. In continuum mechanics, the energy equation for an inviscid gas can be written as [12]: ${{\partial \,s }/{\partial \xi ^1}}\,=0$, where  $s$, where s is entropy. In this case,  
$\psi \,=\,s$, $A_1\,=\,0$. The corresponding equation for the law of conservation of momentum is also given in [12].]}

Equations (6) and (7) can be collapsed into the relation
$$
d\psi\,=\,A_{\mu }\,d\xi ^{\mu }\eqno(8)
$$ 
(Here the summation is carried out over the repeating index.)
This relationship can be written in the following form:
 
This relation can be written in the following form:
$$
d\psi \,=\,\omega\eqno(9) 
$$
where  $\omega \,=\,A_{\mu }\,d\xi ^{\mu }$ is a skew-symmetric differential form of the first degree [2]. Since the equations of conservation laws are evolutionary, the resulting relationship will also be evolutionary. 

If we also take into account the laws of conservation of angular momentum and mass, then the evolutionary relation will take the form 
$$
d\psi \,=\,\omega^p\eqno(10)
$$ 
where the degree of the form $p$ takes the values  $p\,=\,1,2,3$.   

\subsection*{Properties of the evolutionary relation.} 

Here you should pay attention to the fact that the resulting relation is expressed through skew-symmetric forms [1-3]. It is the formalism of skew-symmetric forms that makes it possible to reveal the features of differential mathematical physics equations. 

In the works [1-3], as well as in Appendix 1, it is shown that in addition to external skew-symmetric forms that are defined on integrable manifolds or on integrable structures, there are skew-symmetric forms whose basis are non-integrable manifolds [13]. As it is known, closed external forms are differentials and, therefore, have invariant properties. And skew-symmetric forms defined on non-integrable manifolds are evolutionary forms and have a unique property. They can generate closed external forms that have invariant properties and describe invariant structures. 

Such properties of skew-symmetric differential forms make it possible to reveal the properties of the evolutionary relation and its role in the mathematical physics equations. 

\bigskip 

As was shown, evolutionary relation (9) was obtained on the accompanying manifold (the manifold formed by the trajectories of the described environment), which is a deformable non-integrable manifold. (Examples of non-integrable spaces of differential equations are tangent spaces of differential equations, accompanying manifolds of differential equations, Lagrangian manifolds.) 

In relation (9), the form $\omega \,=\,A_{\mu }\,d\xi ^{\mu }$ is an evolutionary form, since it is defined on a non-integrable manifold. 

The differential of the evolutionary form (which is defined on a non-integrable manifold) is not equal to zero, since it contains a non-zero differential of the metric form of the non-integrable manifold [13]. This means that the  evolutionary form is not closed skew-symmetric form, and, therefore, is not a differential. 

It turns out that on the left in the evolutionary relation (9) there is a differential, and on the right is an evolutionary form that is not a differential. This indicates that evolutionary relation (9) on the original accompanying manifold turns out to be non-identical. 

Evolutionary relation (10) on the original manifold is also a non-identity relation. 

\bigskip 

The non-identical evolutionary relation has one more property. It is a self-changing relation, since when one term changes, the second term cannot be compared with it and the process continues. 

The evolutionary relation, which is obtained by studying the consistency of the equations of conservation laws, allows us to study the integrability of the mathematical physics equations. 

\subsection{Nonintegrability of mathematical physics equations on the original tangent space.} 

A non-identical evolutionary relation cannot be directly integrated, since on the right there is an evolutionary form that is not a closed form, that is, a differential. (The evolutionary form is not a closed form, since the commutator of the evolutionary form, through which the differential is expressed, includes a non-zero commutator of the metric form of the non-integrable manifold.)

The mathematical physics equations cannot be collapsed into an identity relation and directly integrated. This indicates that on tangent space the mathematical physics equations turn out to be non-integrable. (Here we can note the following. Since physical processes in the described media or systems occur as a result of various non-potential influences, for example, such as energy, force, electromagnetic and so on, the mathematical physics equations must contain non-potential terms. This is precisely the reason for the non-integrability of the mathematical physics equations on the original tangent space.) 

The non-integrability of the mathematical physics equations means that the solutions to the mathematical physics equations in this case cannot be analytical solutions, that is, functions. They will depend on the non-zero differential of the evolutionary form, that is, on the commutator of the evolutionary form, and not just on the variables. The derivatives of such solutions do not form a differential. 

Such solutions, which can be obtained by numerical methods, have a physical meaning. They describe the nonequilibrium state of the environment. 

[The evolutionary relation includes the functional  $\psi$, which characterizes the co-state of the described environment. The presence of a state functional differential $d\psi$ in the evolutionary relation indicates the presence of a state function, which corresponds to the equilibrium state of the environment. But since the evolutionary relation is non-identical, it is impossible to obtain the differential of the state functional 
$d\psi$ from it. This indicates the absence of a state function and means that the state of the described environment is nonequilibrium. (It is obvious that the internal forces that cause disequilibrium must be described by a commutator of an open evolutionary form $\omega$). It turns out that the non-integrability of the equation of mathematical physics indicates a nonequilibrium state of the medium.
Another property of a non-identical evolutionary relation, namely its self-change, indicates that the non-equilibrium state of the material environment turns out to be self-changing. The state of the material environment changes, but remains nonequilibrium, since the evolutionary relationship remains non-identical in the process of self-change.] 

It should be noted here that the evolutionary form and evolutionary relation were obtained from the mathematical physics equations when studying the consistency of conservation law equations. The non-closedness of the evolutionary form and the non-identity of the evolutionary relationship mean that the equations of conservation laws turn out to be inconsistent. It is the inconsistency of the equations of conservation laws that is the reason for the non-integrability of the mathematical physics equations on the original tangent space. (As was shown in Chapter 1, the non-integrability of differential equations is also associated with the inconsistency of the derivatives of the differential equation.) 

\subsection{Realization of integrability of mathematical physics equations.}

The mathematical physics equations can become integrable if the evolutionary form, which is not a differential, produces a closed external form, which is a differential. In this case, the non-identical evolutionary relation will yield an identical relation with a closed external form (that is, a differential), which can be directly integrated. This will indicate the integrability of the equation of mathematical physics. 

But since the differential of a closed external form is equal to zero, and the differential of an evolutionary form is not equal to zero, then the transition from an evolutionary form to a closed external form can only occur under a degenerate transformation, that is, under a transformation that does not preserve the differential.

\subsection*{Degenerate transformation. Integrability of mathematical physics equations on a cotangent manifold.}

A degenerate transformation is a transition from the differential of the evolutionary form $\omega$, which is not equal to zero, to the differentials of the closed inexact external form and the corresponding dual form, which are equal to zero. This can be written as:

$d\omega\ne 0 \to $ (вырожденное преобразование) $\to d_\pi \omega=0$,
$d_\pi{}^*\omega=0$ 

Here the condition $d_\pi{}^*\omega=0$, which is the condition for the dual form to be closed, indicates the realization of the closed dual form  ${}^*\omega_\pi$, describing some integrable structure, that is, a pseudostructure $\pi$.  

The implementation of the dual form, describing the integrable structure, indicates the emergence of an integrable structure. 

And condition $d_\pi \omega=0$, which is the condition for the closedness of the imprecise external form, indicates the implementation on the integrable structure of the closed inexact external form $\omega_\pi$ , which is a differential, that is, a conserved quantity.

That is, with a degenerate transformation from an evolutionary form, which is defined on a non-integrable manifold, a closed inexact external form is obtained, which is defined on an integrable structure. In this case, a transition occurs from a non-integrable manifold with an evolutionary form to a pseudostructure (which is described by a closed dual form) with a closed inexact external form. (Such integrable structures are structures of a cotangent manifold.)

{\footnotesize Mathematically, a degenerate transformation is realized as a transition from one coordinate system to another non-equivalent coordinate system: from an accompanying non-inertial coordinate system (associated with a manifold, formed by the trajectories of elements of the described medium) to a locally inertial one, defined on integrable structures or manifolds. An example is the transition from the Lagrangian coordinate system to the Eulerian one.]} 

\bigskip 
As a result of a degenerate transformation on an integrable structure, that is, on a pseudostructure, from the evolutionary relation (9) we obtain the relation  
$$
d_\pi\psi=\omega_\pi\eqno(11)
$$
which turns out to be identical since it contains only differentials (the closed external form $\omega_\pi$ is an internal differential).  

\bigskip
Since the identity relation consists of differentials, it can be integrated. This indicates that the equation of mathematical physics becomes integrable. Solutions to mathematical physics equations on integrable structures (pseudostructures) are functions. But such analytical solutions are discrete functions, since they are defined only on structures. These are the so-called generalized solutions. 

Here it should be emphasized that the mathematical physics equations, which are a system of equations, as well as in the case of a single equation (as shown in Chapter 1), have double solutions: on the original coordinate space and on the integrable structures. 

Moreover, as shown, they are defined on different spatial objects. Such solutions cannot be obtained by continuous simulation. 

As was said, solutions that are not functions can only be obtained by numerical methods. And solutions of the second type, which are discrete functions, can be obtained by analytical methods or numerical methods. But unlike solutions of the first type, which are not functions and are obtained by numerical methods in a non-inertial coordinate system on the original coordinate space, generalized solutions are found using a locally inertial coordinate system defined on integrable structures. It turns out that to obtain double solutions of the mathematical physics equations by numerical methods it is necessary to use two nonequivalent coordinate systems. 

It turns out that the degenerate transformation brings the equation of mathematical physics to an integrable form and to the implementation of discrete solutions.

{\footnotesize [As was shown above, the non-integrability of the mathematical physics equations corresponds to a nonequilibrium state of the medium. This is due to the fact that from a non-identical evolutionary relation it is impossible to obtain a differential state functional, the presence of which would indicate an equilibrium state of the environment. 

When realizing the integrability of an equation of mathematical physics, from a non-identical evolutionary relation we obtain an identical relation, from which we can obtain the differential of the state functional d?. This indicates the presence of a state function, and the transition of the material environment from a nonequilibrium state to a locally equilibrium state. (At the same time, the general state of the material system remains nonequilibrium.) 

The transition from a nonequilibrium state to a locally equilibrium state means that the immeasurable quantity, which was described by the commutator of the evolutionary form and acted as an internal force, turns into a measurable, potential quantity of the material environment. This is associated with the emergence of physical structures and is manifested in the appearance of some observable formations in the material environment. Waves, vortices, oscillations, turbulent pulsations, etc. are such formations [14]. This is described by the solution on the realized integrable structure.]}  

Thus, it was shown that the integrability of mathematical physics equations, on which no additional conditions are imposed, is realized in the process of solving mathematical physics equations. {\bf And this happens using a degenerate transformation,} which carries out the transition from a non-integrable space to an integrable manifold and is accompanied by the emergence of an integrable structure. 

\bigskip

A degenerate transformation can be realized only under additional conditions. Such additional conditions are the vanishing of such functional expressions as Jacobians, determinants, Poisson brackets, residues, etc. Such conditions can be realized if there are any degrees of freedom. This happens discretely in the process of solving the mathematical physics equations, which is described by a self-changing evolutionary relation. 

These conditions describe an integrable structure to which a closed dual form corresponds. 

The realization of these conditions leads to the fact that from the tangent space on which the evolutionary form is defined, an integrable structure, that is, a closed dual form, is realized. This leads to the fact that on an integrable structure, from an evolutionary form a closed external form is obtained, conjugate to the dual form. 

In this case, a transition occurs from a tangent non-integrable space to an integrable cotangent manifold, which is formed by integrable structures. 

The peculiarity of a degenerate transformation is that the integrable structure cannot be the structure of a tangent non-integrable manifold. From such integrable structures a cotangent manifold is formed. (Such structures are sections of a cotangent manifold.) 

\bigskip

Degenerate transformations not only bring the mathematical physics equations to an integrable form, but they reveal the quantum properties of the mathematical physics equations and the mechanism of the emergence of quantum structures.

\section{Quantum properties of mathematical physics equations. Generation of quantum structures} 

It was shown that when implementing degenerate transformations, in which the mathematical physics equations become integrable, integrable structures (described by dual forms) with differentials (closed inexact external forms that are a conserved under non-degenerate transformation quantity) arise. 

And here there is a specific feature. 
When implementing a degenerate transformation, only a mini structure can be formed. Integrable structures are formed by such mini-structures that arise with each successive implementation of degenerate transformations. 

This reveals the quantum properties of the mathematical physics equations. 

What physical meaning do such integrable structures with conserved quantities carry? 

As it has been shown, with a degenerate transformation, a dual form is realized that describes the integrable structure, and a closed inexact external form, which is a differential, that is, a conserved, invariant quantity. 

From the conjugacy of a closed imprecise external form and a dual form it follows that a closed imprecise external form and the corresponding dual form a differential geometric structure that describes a physical structure: an integrable structure with a conserved (under non-degenerate transformation) quantum. 

Here, the conserved quantity is the quantity that is conserved in the resulting quantum. And the other quantum that arises will have its own conserved quantity.

It should be noted that the value of the conserved quantity is obtained by solving the equations of mathematical physics.

Since integrable structures are formed by mini-structures, this indicates that the physical structure consists of quanta: mini-structures with a conserved (under non-degenerate transformation) quantity that arise with each successive implementation of a degenerate transformation. 

{\footnotesize [Here it should be noted that the generator of emerging quanta are solutions of non-integrable mathematical physics equations. As shown, the emergence of integrable structures occurs with the help of a degenerate transformation during the transition from the tangent non-integrable space of mathematical physics equations to the cotangent, integrable manifold. This indicates that the generator of emerging quanta are solutions of non-integrable equations, which are not functions, and can only be obtained by numerical methods in the original coordinate space.]}

Since the physical structure contains a differential (closed external form), it turns out to be invariant under all transformations preserving the differential (that is, under non-degenerate transformations). 

It turns out that physical structures: emerging invariant structures with conserved quantity describe invariant objects. At the same time, the implementation of mini-structures with conservation quantities (quanta) indicates the quantum nature of physical structures and reveals the quantum structure of invariant objects.

{\footnotesize [Above we considered a degenerate transformation with skew-symmetric forms of the first degree, which is obtained by the interaction of two equations of conservation laws (the equation of conservation of energy and the equation of conservation of momentum). When other equations of conservation laws are taken into account, skew-symmetric forms of degree $p$ are obtained (see (10)). In this case, degenerate transformations with skew-symmetric forms of degree $p$ are realized, under which the resulting invariant structures will be surfaces of the corresponding dimension. These can be eikonals, potential surfaces, etc. This indicates that the mathematical physics equations can describe the emergence of various invariant objects. These are spatial objects that carry information.]} 

Invariant objects appear in many invariant mathematical formalisms. Such objects are potential, quantum objects, physical structures that form physical fields [14] and so on. The quantum nature of the emergence of physical structures indicates the quantum structure of invariant objects. It can be noted that the emergence of invariant objects occurs through the interaction of two or more quantities of different natures. This is explained by the fact that they are obtained by interacting two or more equations.

Quantum properties of the mathematical physics equations allow the mathematical physics equations to describe the processes of the emergence of various structures and discrete transitions.

\bigskip
{\Large\bf Conclusion}

The work shows that the mathematical physics equations that describe physical processes, and on which no additional conditions are imposed, have quantum properties. This is due to the integrability of differential equations, which depends on the consistency of the derivatives and the consistency of the equations if the mathematical physics equations are a system of equations. 

A study of the consistency of derivatives and mathematical physics equations has shown that the implementation of integrability of differential equations is accompanied by the emergence of various integrable structures, examples of which are characteristics, eikonals, potential surfaces, etc. This is described by a degenerate transformation that does not preserve the differential. An example of such a transformation is the Legendre transformation.

A degenerate transformation has a unique property. 

During a degenerate transformation, integrable structures arise and a transition occurs from a non-integrable space to the resulting structures of an integrable manifold. This is quantum in nature. With each implementation of a degenerate transformation, mini structures arise, from which integrable structures are formed with successive implementations of degenerate transformations. 

Degenerate transformations occur in the process of solving mathematical physics equations when implementing any degrees of freedom. This corresponds to the vanishing of such functional expressions as determinates, Jacobians, Poisson brackets, residues, etc. 

The process of implementing degenerate transformations of the mathematical physics equations reveals the quantum nature of the mathematical physics equations, which makes it possible to describe discrete transitions, the emergence of various structures and invariant objects. 

It can be assumed that the basis of quantum processes is a mechanism described by the properties of degenerate transformations of skew-symmetric forms: the emergence of quantum structures when implementing any degrees of freedom of the medium or physical system under study and the transition from non-integrable spaces with non-potential quantities to the structures of an integrable manifold.

\bigskip
Such properties of the mathematical physics equations were obtained using skew-symmetric differential forms. Moreover, in addition to external skew-symmetric forms, the basis of which are integrable manifolds and structures, evolutionary forms were used on non-integrable manifolds, which have non-traditional operators: non-identity relations and degenerate transformations. 

\bigskip
{\footnotesize
\rightline{\large\bf Appendix 1}

{\large\bf Some properties of skew-symmetric forms.} 

\bigskip

{\bf Closed inexact exterior forms.}

The external differential form of the degree $p$ ($p$-form) can be written
as [2]:
$$
\theta^p=\sum_{i_1\dots i_p}a_{i_1\dots i_p}dx^{i_1}\wedge
dx^{i_2}\wedge\dots \wedge dx^{i_p}\quad 0\leq p\leq n\eqno(1)
$$
Here $a_{i_1\dots i_p}$ is the function of independent variables $x^1, ..., x^n$,
$n$ is the space dimension, and
$dx^i$, $dx^{i}\wedge dx^{j}$, $dx^{i}\wedge dx^{j}\wedge dx^{k}$, \dots\
is the local basis subject to the condition of skew-symmetry:
$$
\begin{array}{l}
dx^{i}\wedge dx^{i}=0\\
dx^{i}\wedge dx^{j}=-dx^{j}\wedge dx^{i}\quad i\ne j
\end{array}\eqno(2)
$$
The exterior form differential $\theta^p$ is expressed by the formula
$$
d\theta^p=\sum_{i_1\dots i_p}da_{i_1\dots
i_p}\wedge dx^{i_1}\wedge dx^{i_2}\dots \wedge dx^{i_p} \eqno(3)
$$
The form called as a closed one if its differential equals to zero:
$$
d\theta^p=0\eqno(4)
$$
From condition (4) one can see that the closed form
is a conservative quantity.
This means that such a form can correspond to the conservation law
(for physical fields), i.e. a conservative quantity.

If the form be closed only on pseudostructure, i.e. this form is a closed
inexact one, the closure condition can be
written as
$$
d_\pi\theta^p=0\eqno(5)
$$
In this case the pseudostructure $\pi$ obeys the condition
$$
d_\pi{}^*\theta^p=0\eqno(6)
$$                    	    				
here ${}^*\theta^p$ is the dual form. 

\bigskip
{\bf Features of closed exterior forms.} 
 
The closed outer form is a differential [1] and therefore an invariant in all transformations that conserve the differential (i.e., in degenerate transformations). In this case, the closed inaccurate exterior form is an interior (on the pseudo structure) differential.
From the closure condition it follows that the closed form is a conserved (under non-degenerate transformations) value and satisfies the exact conservation law (as a conserved variable).

\bigskip

{\bf Evolutionary skew-symmetric differential forms.}

Evolutionary skew-symmetric differential forms are skew-symmetric forms, based on non-integrable manifolds.
The evolutionary form can be expressed similarly to the exterior differential form. However, a member associated with the differential of the evolutionary form basis appears in the differential of the evolutionary form, since such a form is based on a non-integrable manifold [2].

The evolutionary form differential takes the form
$$
d\theta^p{=}\!\sum_{i_1\dots\i_p}\!da_{i_1\dots\i_p}\wedge
dx^{i_1}\wedge dx^{i_2}\dots \wedge
dx^{i_p}{+}\!\sum_{i_1\dots\i_p}\!a_{i_1\dots\i_p}
d(dx^{i_1}\wedge dx^{i_2}\dots \wedge dx^{i_p})\eqno(7)
$$
where the second term is connected with the basis differential being
nonzero:
$d(dx^{i_1}\wedge dx^{i_2}\wedge \dots \wedge dx^{i_p})\neq 0$.
(For the exterior form defined on integrable manifold one has
$d(dx^{i_1}\wedge dx^{i_2}\wedge \dots \wedge dx^{i_p})=0$).
[Then, the exterior multiplication symbol (?) and the summation sign ? will be omitted and it will be assumed that summation is done with the double index.]

\bigskip
{\bf Features of evolutionary forms.}

Since the differential of an evolutionary form is not equal to zero, then the evolutionary form cannot be a closed external form, that is, {\bf cannot be a differential}.

But from the evolutionary form closed inaccurate external forms can be obtained, which are differentials. And this happens using a degenerate transformation, that is, using a transformation that does not preserve the differential.

\bigskip
{\bf Degenerate transformations and non-degenerate transformations of skew-symmetric differential forms.}

Degenerate transformations, and non-degenerate transformations (transformations that conserve the differential) are the transformations of the theory of skew-symmetric differential forms. 
At the same time, the non-degenerate transformations are the transformations of exterior skew-symmetric differential forms [2], the basis of which are integrable manifolds, and the degenerate transformations, except for exterior skew-symmetric forms, are also associated with skew-symmetric forms, which are evolutionary forms [3], as they were defined on non-integrable manifolds. 

The peculiarity of non-degenerate transformations is that non-degenerate transformations
carry out transitions between integrable structures on which differentials (closed external forms) are defined. And this happens on integrable manifolds.

And with a degenerate transformation, an integrable structure arises and a transition occurs from a non-integrable manifold with an evolutionary form, the differential of which is not equal to zero, to the integrable structure of an integrable manifold with a closed inexact external form, the differential of which is equal to zero. In this case, the integrable structure is described by a dual form.

When implementing a degenerate transformation from an evolutionary skew-symmetric form, the basis of which is a non-integrable manifold, a closed inexact external form is realized on an integrable structure. Such a process, which is accompanied by the emergence of an integrable structure, has a quantum character.]}

[1] {\it Cartan E.} Les Systems Differentials Exterieus ef Leurs Application
Geometriques. -Paris, Hermann, 1945.

[2] {\it Petrova L.I.} Exterior and evolutionary differential forms in mathematical physics: Theory and Applications, -Lulu.com,157, 2008. 

[3] {\it Petrova L.I.} Role of skew-symmetric differential forms in mathematics, 2010.  
http://arxiv.org/pdf/1007.4757.pdf. 

[4] {\it Petrova L.} Qualitative Investigation of Hamiltonian Systems by Application of Skew-Symmetric Differential Forms, //Symmetry,  MDPI (Basel, Switzerland), 2021, Vol 13,  No.1, 1-7.

[5] {\it Р. Курант и Д. Гильберт} Методы математической физики, М.: ГИТТЛ, том 2, 1951. 

[6] {\it Clark J F., Machesney M.} The Dynamics of Real Gases, Butterworths, London, 1964.  

[7] {\it В.И. Смирнов} Тех-Теор. Лит., Москва, V.4., 1957.

[8] {\it R. Courant} Partial Differential Eguations, New York.London, 1962.

[9] {\it W. Pauli} Theory of Relatiyity, Pergamon Press, 1958. 

[10] {\it А.Н. Тихонов, А.А. Самарский} Уравнения математической физики, Тех-Теор. Лит., Москва, 1953.

[11] {\it Petrova L.} Connection between functionals of the field-theory equations and state functionals 
of the mathematical physics equations. //Journal of Physics: Conference Series, IOP Publishing ([Bristol, UK], England), 2018,  Vol 1051, No. 012025, 1- 8. 

[12] {\it Petrova L.I.} Hidden properties of the Navier-Stokes equations. Double solutions. Origination of turbulence. //Theoretical Mathematics and Applications (TMA), 2014, Vol 4, No. 3, 91-108.

[13] {\it Tonnelat M.-A.} Les principles de la theorie electromagnetique et la relativite, Masson, Paris, 1959.

[14] {\it Petrova L.I.} Hidden Unique Possibilities of Mathematical Physics Equations (The Formalism of Skew-Symmetric Forms). // Computational Mathematics and Modeling, 2022, 33, pp. 121-135.

\end{document}